\documentclass[12pt]{amsart}
\usepackage{amssymb}
\usepackage[PostScript=dvips]{diagrams}

\newtheorem{thm}{Theorem}[section]
\newtheorem{lem}[thm]{Lemma}
\newtheorem{cor}[thm]{Corollary}
\newtheorem{prop}[thm]{Proposition}

\newtheorem{dfn}[thm]{Definition}

\newtheorem{rem}[thm]{Remark}

\def\square{\vbox{
      \hrule height 0.4pt
      \hbox{\vrule width 0.4pt height 5.5pt \kern 5.5pt \vrule width 0.4pt}
      \hrule height 0.4pt}}

\def\id{\mathop{\rm id}\nolimits}

\def\Ker{\mathop{\rm Ker}\nolimits}

\def\min{\mathop{\rm min}\nolimits}

\def\max{\mathop{\rm max}\nolimits}

\newcommand{\Z}{{\mathbb Z}}

\newcommand{\R}{\ensuremath{\mathbb R}}
\newcommand{\Q}{\ensuremath{\mathbb Q}}

\newcommand{\calB}{\ensuremath{\mathcal{B}}}

\newcommand{\calZ}{\ensuremath{\mathcal{Z}}}

\newcommand{\calR}{\ensuremath{\mathcal{R}}}

\let\la=\langle
\let\ra=\rangle

\begin{document}

\newcommand{\auths}[1]{\textrm{#1},}
\newcommand{\artTitle}[1]{\textsl{#1},}
\newcommand{\jTitle}[1]{\textrm{#1}}
\newcommand{\Vol}[1]{\textbf{#1}}
\newcommand{\Year}[1]{\textrm{(#1)}}
\newcommand{\Pages}[1]{\textrm{#1}}

\address{Department of Mathematics\\
National University of Singapore\\
Singapore 117543\\
Republic of Singapore\\
matwuj@nus.edu.sg} 
\author{ Jie Wu}
\title{A Braided Simplicial Group}
\thanks{Research is supported in part by the Academic Research Fund of the National University of Singapore RP3992646}
\begin{abstract}
By studying braid group actions on Milnor's construction of the $1$-sphere, we show that the general higher homotopy group of the $3$-sphere is the fixed set of  the pure braid group action on certain combinatorially described group. We also give certain representation of higher differentials in the Adams spectral sequence for $\pi_*(S^2)$.
\end{abstract}

\maketitle

\section{Introduction}
In this article, we study the homotopy groups by considering the braid group actions on simplicial groups. The point of view here is to establish a relation between the fixed set of braid group actions and the homotopy groups of the $3$-sphere. We first recall a combinatorial description of the homotopy groups of the $3$-sphere in~\cite{Wu1}.

Let $F(x_1,\cdots,x_n)$ be the free group generated by the letters
$x_1,\cdots,x_n$. Let $w(x_1,\cdots,x_n)=x_{i_1}^{\epsilon_1}\cdots
x_{i_t}^{\epsilon_t}$ be a word. Given $a_1,\cdots,a_n\in 
F(x_1,\cdots,x_n)$, we
write $w(a_1,\cdots,a_n)\in F(x_1,\cdots,x_n)$ for 
$a_{i_1}^{\epsilon_1}\cdots
a_{i_t}^{\epsilon_t}$. The {\it $n$-th $W$-group $G(n)$} is the quotient 
group of
$F(x_1,\cdots,x_n)$ modulo the following relations:

\begin{enumerate}
\item[$(\calR_1)$] the  product $x_1\cdots x_n$;
\item[$(\calR_2)$] the words $w(x_1,\cdots,x_n)$ that satisfy:
$
w(x_1,\cdots,x_{i-1},1,x_{i+1},\cdots,x_n)=1
$
for $1\leq i\leq n$.
\end{enumerate}

Relations $\calR_2$ consist of all of words that will be trivial if one 
of the
generators is replaced by the identity element $1$. The smallest
normal subgroup of $F(x_1,\cdots,x_n)$ which contains relations 
$\calR_1$ and $\calR_2$ was determined 
 as a subgroup of $F(x_1,\cdots,x_n)$ 
generated by
certain systematic and uniform iterated commutators~\cite{Wu1}.

\begin{thm}~\cite{Wu1}\label{theorem1.1}
For $n\geq 3$, the homotopy group $\pi_n(S^3)$ is isomorphic to the 
center of
$G(n)$.
\end{thm}

A natural question arisen from Theorem~\ref{theorem1.1} is how to give a group theoretical approach to the homotopy groups, that is how to understand the center of the group $G(n)$. There is a canonical braid group action on $G(n)$ which is induced by the canonical braid
group action on the free group $F(x_1,\cdots,x_n)$, namely
$$
\sigma_i(x_j)=\left\{
\begin{diagram}
x_{i+1}&if&j=i&\\
x_{i+1}^{-1}x_ix_{i+1}&if&j=i+1&\\
x_j&&otherwise\\
\end{diagram}
\right.
$$
for $1\leq i\leq n-1$.  These actions gives a canonical homomorphism from
the braid group $B_n$ into the automorphism group of $G(n)$.
Since Quillen's  plus construction of the classifying space for
the stable braid group is  (up to homotopy type)
the double loop space of the 3-sphere~\cite{Cohen}, these braid group actions do not seem occasional event. Fred Cohen therefore
conjectured that the center of $G(n)$ is the  fixed set of the
braid group action on $G(n)$.  We answer Cohen's question as follows.
Let $K_n$ be the pure braid group, that is, $K_n$ is the normal divisor of
$B_n$ generated by $\sigma_1^2$ (See~\cite{MKS}). Equivalently $K_n$ is the kernel of the canonical homomorphism from $B_n$ to the symmetric group $\Sigma_n$. In geometry, the group $K_n$ is the fundamental group of the configuration space $F(\R^2,n)$, where 
$$
F(M,n)=\{(x_1,\cdots,x_n)\in M^n|x_i\not=x_j\,\, {\rm for}\,\, i\not=j\}
$$
for any manifold $M$ (See~\cite{Cohen}).  
Let $Z(G(n))$ be the center of 
$G(n)$.

\begin{thm}\label{theorem1.2}
For $n\geq4$, then
\begin{enumerate}
\item[1)] the center of $G(n)$ is the fixed set of the pure braid group
action on $G(n)$ and so is $\pi_n(S^3)$;
\item[2)] the fixed set of the braid group action on $G(n)$ is the subgroup
$$
\{ x\in Z(G(n))|2x=0\}.
$$
\end{enumerate}
\end{thm}
 
We should point out that the determination of the fixed set of $K_n$-action on $G(n)$ by (combinatorial) group theoretic means seems beyond the reach of current techniques. On the other hand, braid group actions have largely studied in several areas such as group theory and low dimensional topology. Various problems arising from physics are related to braid group actions as well. Theorem~\ref{theorem1.2} suggests that the homotopy groups play certain role for braid group actions. In the range in which $\pi_*(S^3)$ is known ($*\leq55$, see~\cite{Mahowald,T}),
by homotopy theoretic means, we gain insight  into these difficult group 
theoretic questions. 

The article is organized as follows. In section~\ref{section2}, we study the
braid group action on the Milnor's construction on the simplicial $1$-sphere.
A relation between the simplicial structure and the braid group action is given
in Proposition~\ref{proposition2.1}. This relation is based on the direct calculation. Roughly the braid group action interchanges the faces together with conjugates. It is possible to have a more general theory to study group actions on simplicial groups, particularly simplicial group models for iterated loop spaces. But we only intend to investigate the most important example $F(S^1)$ in this article. After establishing the systematic relations between the braid group actions and the simplicial structure, braided simplicial
groups are introduced in this section. Then we show that loop simplicial group and 
the Moore-Postnikov system of a braided simplicial group are braided. 
Theorem~\ref{theorem2.8} and Proposition~\ref{proposition2.9} give a relation 
between the fixed set of the braid group action and the homotopy groups
for a general braided simplicial group. Theorem~\ref{theorem1.2} follows from 
Lemma~\ref{lemma2.7} and Theorem~\ref{theorem2.11}. In section~\ref{section3},
we study a braided representation of the Milnor's construction of the 
simplicial $1$-sphere into a simplicial algebra. Theorem~\ref{theorem3.8}
gives certain representation of  higher differentials in the Adams spectral
sequence for $\pi_*(S^2)$.  

The author would like to thanks Professors Fred Cohen and Jon Berrick for their helpful suggestions and kind encouragements.

\section{Braid Group Actions on $F(S^1)$}\label{section2}
\subsection{Braided simplicial groups}
Let $F(S^1)$ be Milnor's construction of the simplicial $1$-sphere $S^1$.
Then $F(S^1)_{n+1}=F(y_0,\cdots,y_n)$ is a free group generated by
letters $y_0,\cdots,y_n$ with the 
following simplicial structure
$$
d_jy_k=\left\{\begin{array}{cccccc}
y_{k-1}&\mbox{ for}& \mbox{ $j\leq k,$}\\
1&\mbox{ for}&\mbox{ $j=k+1,$}\\
y_k&\mbox{ for} &\mbox{ $j>k+1,$}\\
\end{array}
\right.
$$
and
$$
s_jy_k=\left\{\begin{array}{cccccc}
y_{k+1}&\mbox{ for}&\mbox{ $j\leq k,$}\\
y_ky_{k+1}&\mbox{ for}&\mbox{ $j=k+1,$} \\
y_k&\mbox{ for}&\mbox{ $j>k+1$,}\\
\end{array}
\right.
$$
for $0\leq j\leq n+1$, where $y_{-1}=(y_0\cdots y_{n-1})^{-1}$ in $F(S^1)_n$ 
(See~\cite[Lemma 4.1]{Wu1}). Let the braid group $B_{n+1}$ act on 
$F(S^1)_{n+1}=F(y_0,\cdots,y_n)$ in the usual way, that is,
$$
\sigma_k(y_j)=\left\{
\begin{array}{ccc}
y_{k+1}&\quad {\rm if}&\quad j=k\\
y_{k+1}^{-1}y_ky_{k+1}&\quad {\rm if}&\quad j=k+1\\
y_j & &{\rm otherwise}\\
\end{array}
\right.
$$
for $0\leq k\leq n-1$. Let $\sigma_{-1}$ be an automorphism of 
$F(y_0,\cdots,y_n)$ defined by
$$
\sigma_{-1}(y_0)=y_0^{-1}y_{-1}y_0=y_0^{-1}y_n^{-1}\cdots y_1^{-1}\quad
\sigma(y_j)=y_j\quad{\rm for}\quad j>0.
$$
The subgroup of the automorphism group of $F(y_0,\cdots,y_n)$ generated by
$\sigma_j$ for $-1\leq j\leq n-1$ is the braid group $B_{n+2}$. 
By direct calculation, we have
\begin{prop}\label{proposition2.1}
The following identities hold for the braid groups action on $F(S^1)$:
\begin{equation}\label{equation1}
d_j\sigma_k=\left\{
\begin{diagram}
\sigma_{k-1}d_j &\quad j\leq k\\
d_{k+2} &\quad j=k+1\\
d_{k+1} &\quad j=k+2\\
\sigma_k d_j &\quad j>k+2\\
\end{diagram}
\right.
\end{equation}
\begin{equation}\label{equation2}
s_j\sigma_k=\left\{
\begin{diagram}
\sigma_{k+1}s_j&\quad  j\leq k\\
\sigma_{k+1}\circ \sigma_{k}\circ s_{k+2}&\quad j=k+1\\
\sigma_k\circ\sigma_{k+1}\circ s_{k+1}&\quad j=k+2\\
\sigma_k s_j&\quad j>k+2\\
\end{diagram}
\right.
\end{equation}
\end{prop}

By this proposition, we give the following definition.
\begin{dfn}\label{definition2.2}{\rm
A simplicial group $G$ is called {\it braided} if there is a braid group 
$B_{n+1}$-action on $G_n$ for each $n$ such that the 
identities~\ref{equation1} and~\ref{equation2} are satisfied, where $B_{n+1}$ 
is considered as the braid group generated by $\sigma_{-1},\sigma_0,
\cdots, \sigma_{n-2}$.}
\end{dfn}

Let $G$ be a simplicial group and let
 $NG$ be the Moore chain complex of $G$, that is,  
$$
NG_n=\{x\in G_n|d_jx=1\quad {\rm for}\quad j>0\}
$$. Let $\calZ(G)$ and $\calB(G)$ be the 
sets of cycles and boundaries of $G$, respectively, that is,
$$
\calZ_n(G)=\{x\in G_n|d_jx=1\quad {\rm for\quad  all}\quad j\},
$$
$$
\calB_n(G)=\{d_0x|x\in N_{n+1}(G)\}.
$$
By Moore's classical theorem~\cite{Moore}, $\pi_n(G)=\calZ_n(G)/\calB_n(G)$. Let $G$ be
a braided simplicial group.
 A subgroup $H$ of $G$ is called a
{\it braided subgroup} if $\sigma(H)\subseteq H$ for any $\sigma\in B_{n+1}$, 
that is, $H$ is invariant under $B_{n+1}$-action.

\begin{prop}\label{proposition2.2}
Let $G$ be a braided simplicial group. Then the  subgroups $\calZ(G)_{n}
$ and 
$\calB(G)_{n}$ of $G_n$ are 
braided.
\end{prop}
\begin{proof}
By Proposition~\ref{proposition2.1}, $\calZ(G)_{n}$ is a braided subgroup.
Now let $x=d_0y\in\calB(G)_{n}$, where $y\in NG_{n+1}$. By 
Proposition~\ref{proposition2.1}, we have
$$
\sigma_kx=\sigma_kd_0y=d_0\sigma_{k+1}y
$$
for each $k\geq -1$ and
$$
d_j\sigma_{k+1}y=1
$$
for each $j>0$ and $k\geq-1$. Thus $\sigma_kx\in\calB(G)_{n}$ for  each 
$k\geq-1$ and so $\calB(G)_{n}$ is a braided subgroup, which is the assertion. 
\end{proof}

\noindent{\bf Note:} $N_{n}(G)$ is invariant under the subgroup of 
$B_{n+1}$ generated by $\sigma_k$ with $k\geq0$.
 But it is not invariant under 
$\sigma_{-1}$. 

Since there is a relation
$$
d_0\sigma_{-1}=d_1,\quad d_1\sigma_{-1}=d_0,\quad{\rm and}\quad 
d_j\sigma_{-1}=\sigma_{-1}d_j \quad{\rm for}\quad j>1,
$$
we have

\begin{prop}\label{proposition2.13}
Let $G$ be a braided simplicial group. Then
$$
\calZ(G)_n=NG_n\cap\sigma_{-1}(NG_n).
$$
\
for each $n$.
\end{prop}

Now we study the braided group actions on the Postnikov systems of  a braided
simplicial group $G$. 
Let $I=(i_1,\cdots,i_k)$ be a sequence of non-negative integers and let $d_I$ 
denote the composite of face homomorphisms
$$
d_I=d_{i_1}\cdots d_{i_k}.
$$
Given a simplicial group $G$, the simplicial sub groups $R_nG$ and 
$\bar R_n(G)$ are defined as follows:
$$
R_nG_q=\{x\in G_q| d_I(x)=1\quad{\rm for any}\quad I=(i_1,\cdots,i_{q-n})\},
$$
$$
\bar R_n G_q=\{x\in G_q| d_I(x)\in\calB(G)_n\quad{\rm for any}\quad
I=(i_1,\cdots,i_{q-n})\}.
$$
Let $P_nG=G/R_nG$ and $\bar P_nG=G/\bar R_nG$. Then $\{P_nG\}$ is the 
Postnikov system of $G$ (See~\cite{Curtis,Moore}). The quotient homomorphism
$P_nG\to \bar P_nG$ is a homotopy equivalence (See~\cite{Wu1}). The tower
$$
\cdots\to P_nG\to\bar P_nG\to P_{n-1}G\to\cdots
$$
is called a {\it modified Postnikov system} of $G$. One of the important 
properties of the modified Postnikov system is that the short exact sequence
of simplicial groups
$$
K(\pi_n(G),n)\to \bar P_nG\to P_{n-1}G
$$
is a central extension~\cite[Theorem~2.12]{Wu1}. 

By Propositions~\ref{proposition2.1} and~\ref{proposition2.2}, we have

\begin{thm}\label{theorem2.3}
Let $G$ be a braided simplicial group. Then,
for each $n$, the simplicial quotient groups $P_nG$ and $\bar P_nG$
are braided. Thus the modified Postnikov system of $G$ are braided. 
In particular, there is a braided central extension
$$
K(\pi_n(G),n)\to\bar P_n(G)\to P_{n-1}(G).
$$
\end{thm}

Let $G$ be a braided simplicial group. Then $R_0G$ is a braided simplicial 
subgroup of $G$ by the theorem above. Recall that the  loop simplicial group
$\Omega G$ of $G$ is defined by
$\Omega G_n=\Ker (d_0)\cap R_0G_{n+1}$ with $d_j(\Omega G)=d_{j+1}(G)$ and 
$s_j(\Omega G)=s_{j+1}(G)$ (See~\cite{Curtis}). By 
Proposition~\ref{proposition2.1},
$\Ker(d_0)\cap R_0G_{n+1}$ is invariant under the action of $\sigma_j$ for
$0\leq j\leq n-1$. Let $B'_{n+1}$ be the subgroup of $B_{n+2}$ generated by
$\sigma_0,\cdots,\sigma_{n-1}$. Then $B'_{n+1}\cong B_{n+1}$ under the 
canonical isomorphism which sends $\sigma_j$ to $\sigma_{j-1}$ for $0\leq j\leq n-1$.  Thus we obtained the following theorem.

\begin{thm}\label{theorem2.4}
Let $G$ be a braided simplicial group. Then the loop simplicial group 
$\Omega G$ is braided. Thus any iterated loop simplicial groups of $G$ 
are braided.
\end{thm}

\begin{cor}
The loops and the modified Postnikov system of $F(S^1)$ are braided.
\end{cor}

\subsection{Fixed Sets of Braided Actions}

Let $G$ be a simplicial group and let $x,y\in G_n$ be two elements. We 
call that $x$ is {\it homotopic} to $y$, which is denoted by $x\simeq y$, if
$xy^{-1}\in\calB(G)_n$.

\begin{lem}\label{lemma2.7}
Let $G$ be a braided simplicial group and let $x\in\calZ(G)_n$ be a cycle with
$n\geq1$. 
Then
$$
\sigma_k(x)\simeq x^{-1}
$$
for each $k\geq-1$.
\end{lem}
\begin{proof}
For each $k\geq-1$, consider the element $\sigma_{k+1}s_{k+1}x$. By the 
identities~\ref{equation1} and~\ref{equation2}, we have
$$
d_j\sigma_{k+1}s_{k+1}x=\left\{
\begin{diagram}
\sigma_kd_js_{k+1}x=1&\quad{\rm \quad for}\quad&j<k+1\\
\sigma_kx&\quad{\rm for}\quad&j=k+1\\
d_{k+3}s_{k+1}x=1&\quad{\rm for}\quad&j=k+2\\
d_{k+2}s_{k+1}x=x&\quad{\rm for}\quad&j=k+3\\
\sigma_{k+1}d_js_{k+1}x=1&\quad{\rm for}\quad& j>k+3.\\
\end{diagram}
\right.
$$
Thus $\sigma_kx\simeq x^{-1}$, which is the assertion.
\end{proof}

Let $G$ be a braided simplicial group and let $S$ be a subset of $G_n$. Define
$$
B_{n+1}(S)=\{\sigma\in B_{n+1}|\sigma x\simeq x\quad{\rm for\quad 
 all}\quad x\in 
S\}.
$$
Since $\calB(G)$ is invariant under the braided group action, $B_{n+1}(S)$ is
a subgroup of $B_{n+1}$. Let $\tilde B_{n+1}$ be the kernel of the composite
$$
B_{n+1}\rTo\Sigma_{n+1}\rTo^{\rm sign}\Z/2,
$$
that is, $\tilde B_{n+1}$ is the pre-image of the alternating group $A_{n+1}$.

\begin{thm}\label{theorem2.8}
Let $G$ be a braided simplicial group. Then
\begin{enumerate}
\item[1)] $B_{n+1}(\calZ(G)_n)=\tilde B_{n+1}$ or $B_{n+1}$;
\item[2)] $B_{n+1}(\calZ(G)_n)=B_{n+1}$ if and only if $2\cdot \pi_n(G)=0$.
\end{enumerate}
\end{thm}
\begin{proof}
Let $\sigma\in B_{n+1}(\calZ(G)_n)$ and let $x\in\calZ(G)_n$. By 
Lemma~\ref{lemma2.7}, we have
$$
(\sigma_k^{-1}\sigma\sigma_k)(x)\simeq (\sigma_k\sigma)(x^{-1})
\simeq \sigma_k(x^{-1})\simeq x
$$
for each $k\geq-1$. Thus $B_{n+1}(\calZ(G)_n)$ is a normal subgroup of 
$B_{n+1}$. By Lemma~\ref{lemma2.7}, we have 
$$
\sigma_s\sigma_t\in B_{n+1}(\calZ(G)_n)
$$
for any $s,t\geq-1$. It follows that $B_{n+1}(\calZ(G)_n)=\tilde B_{n+1}$ or
$B_{n+1}$, which is the assertion~1.

If $B_{n+1}(\calZ(G)_n)=B_{n+1}$, then by Lemma~\ref{lemma2.7}
$$
x\simeq x^{-1}
$$ for any $x\in\calZ(G)_n$. Thus $2\cdot\pi_n(G)=0$. Conversely, if 
$2\cdot \pi_n(G)=0$, then
$$
\sigma_k(x)\simeq x
$$
for any $x\in\calZ(G)_n$ and $k\geq-1$. Thus $B_{n+1}(\calZ(G)_n)=B_{n+1}$. 
This shows assertion~2.
\end{proof}

Let $G$ be a braided simplicial group and let $H$ be a subgroup of $B_{n+1}$. 
Define
$$
G_n(H)=\{x\in G_n| \sigma(x)\simeq x\quad{\rm for \quad all}\quad \sigma\in 
H\},
$$
that is, $G_n(H)$ is the ``homotopy'' fixed set of $H$. Then $G_n(H)$
is a subgroup of $G_n$. By Theorem~\ref{theorem2.8}, we have that
$$
 G_n(\tilde B_{n+1})\supseteq \calZ(G)_n.
$$

\begin{prop}\label{proposition2.9}
Let $G$ be a braided simplicial group and let $x\in G_n(\tilde B_{n+1})$ with
$n\geq2$. Then
\begin{enumerate}
\item[1)] $d_j(x)=d_{j+2}(x)$ for each $j$;
\item[2)] $\sigma_kd_j(x)=d_{j+1}(x)$ for each $j,k$;
\item[3)] $d_j(x)$ is a fixed point of $\tilde B_{n}$ for each $j$.
\end{enumerate}
In particular, if $d_jx=1$ for some $j$, then $x\in\calZ(G)_n$.
\end{prop}
\begin{proof}
Since $x\in G_n(\tilde B_{n+1})$, we have
$$
\sigma_s\sigma_t(x)\simeq x
$$
for any $s,t\geq -1$. It follows that
$$
\sigma_{-1}x\simeq \sigma_0 x\simeq \sigma_1x\simeq \cdots\simeq\sigma_{n-2}x.
$$
Now for each $-1\leq k\leq n-3$, we have
$$
d_{k+1}x=d_{k+2}\sigma_k x=d_{k+2}\sigma_{k+1} x=d_{k+3}x
$$
and so assertion~1 follows.

Now for each $0\leq s\leq n-2$, we have
$$
d_{s+1}(x)=d_s\sigma_{s-1} (x)=d_s\sigma_s(x)=\sigma_{s-1}d_s(x).
$$
Assertion~2 follows.

For any $s,t\geq -1$, we have
$$
\sigma_s(\sigma_t d_0(x))=\sigma_s(d_0(\sigma_{t+1}x))=d_0(
(\sigma_{t+1}\sigma_{s+1})x)=d_0x.
$$
Thus
$$
\sigma_t\sigma_s(d_0(x))=d_0(x)
$$
for any $s,t\geq-1$ and so $d_0(x)$ is a fixed point of 
$\tilde B_n$. Since
$$
\sigma_t\sigma_s(d_1(x))=\sigma_s((\sigma_{-1}\sigma_t)(d_0(x)))
=\sigma_s(d_0(x)=d_1(x),
$$
$d_1(x)$ is a fixed point of $\tilde B_{n+1}$ and hence 
assertion~3.
\end{proof}

Let $B_n$ be the subgroup of $B_{n+1}$ generated
by $\sigma_j$ with $j\geq 0$. By inspecting the proof, we have

\begin{prop}\label{proposition2.11}
Let $G$ be a braided simplicial group and let $x\in G_n(\tilde B_{n})$ with
$n\geq2$. Then
\begin{enumerate}
\item[1)] $d_j(x)=d_{j+2}(x)$ for each $j\geq1$;
\item[2)] $\sigma_kd_j(x)=d_{j+1}(x)$ for each $j,k\geq1$;
\item[3)] $d_j(x)$ is a fixed point of $\tilde B_{n-1}$ for each $j$;
\item[4)] $d_0(x)$ is a fixed point of $\tilde B_n$.
\end{enumerate}
In particular, if $d_0x=1$ and  $d_jx=1$ for some $j\geq 1$, 
then $x\in\calZ(G)_n$.
\end{prop}

The following lemma is well-known. We give an elementary proof. Let $S$ be  a subset of a group $G$. We write $\la S\ra$ for the subgroup of $G$ generated by $S$.
\begin{lem}\label{lemma2.10}
Let $w\in F(y_0,\cdots, y_n)$ with $n\geq0$. 
Suppose that there is a positive integer $k$
such that $\sigma_j^k(w)=w$ for $0\leq j\leq n-1$. Then $w$ lies in the
subgroup generated by $y_0y_1\cdots y_n$. In addition, if 
$\sigma_{-1}^kw=w$ and $n\geq1$, then $w=1$.
\end{lem}
\begin{proof}
The proof is given by induction on $n$. The assertion is trivial for $n=0$.
 Let $n=1$.
We may assume that $k=2t$ is an even integer. Let $x_0=y_1^{-1}$ and 
let $x_1=y_0y_1$. Then $F(y_0,y_1)=F(x_0,x_1)$. Since $\sigma_0(y_0)=y_1$
and $\sigma_0(y_1)=y_1^{-1}y_0y_1$, we have  $\sigma_0=\chi_{x_1}$ and so
$$
\sigma_0^k=\chi_{x_1}^t=\chi_{x_1^t}.
$$ 
We can write $w$ as a reduced word in $F(y_0,y_1)=F(x_0)*F(x_1)$. Then
$$w=x_0^{n_1}x_1^{l_1}\cdots x_0^{n_s}x_1^{l_s},$$ 
where $n_j\not=0$ for $2\leq j\leq s$ and $l_j\not=0$
for $1\leq j\leq s-1$. Suppose that $w\not\in \la x_1\ra$. There are two cases: 
$n_1\not=0$ or $n_1=0$. If $n_1\not=0$, then $x_1^tw\not=wx_1^t$. 
This contradicts to that $\chi_{x_1^t}(w)=w$. 
Otherwise, $n_1=0$ and $s>1$.  Then $w=x_1^{l_1}x_0^{n_2}w'$ and
 $x_1^{t+l_1}x_0^{n_2}w'\not=x_1^{l_1}x_0^{l_2}w'x_1^t$. 
This contradicts to that $\chi_{x_1^t}(w)=w$. Hence $w\in \la x_1\ra=\la y_0y_1\ra$.

Now suppose that the assertion holds for $n-1$ with $n>1$.
 Since
$$
F(y_0,\cdots,y_n)=F(y_0,\cdots,y_{n-1})*F(y_n)
$$
is a free product, we can write $w$ as a word
$$
w=y^{l_0}_nw_1y^{l_1}_n\cdots w_ty^{l_t}_n,
$$
where $w_j\not=1\in F(y_0,\cdots,y_{n-1)}$ and $l_j\not=0$ for 
$1\leq j\leq t-1$.
Because
$\sigma_j(y_n)=y_n$ for $j<n-1$, we have
$$
\sigma_j^k(w_i)=w_i
$$
for $1\leq i\leq t$ and $0\leq j\leq n-2$. Let $x=y_0y_1\cdots y_{n-1}$.
 By induction, we have
$$
w_i\in \la x\ra
$$
for $1\leq i\leq t$ and so
$$
w\in \la x,y_n\ra.
$$
Let $q\colon F(y_0,y_1,\cdots,y_n)\to F(y_{n-1},y_n)$ be the projection defined
by $q(y_j)=1$ for $j<n-1$ and $q(y_j)=y_j$ for $j\geq n-1$. Then
$$
q\circ\sigma_{n-1}=\sigma_{n-1}\circ q.
$$
Since $\sigma_{n-1}w=w$, we have $\sigma_{n-1}(q(w))=q(w)$ and so
$$
q(w)\in \la y_{n-1}\cdot y_n\ra.
$$
Because the restriction
$$
q|_{\la x,y_n\ra}\colon \la x,y_n\ra\to F(y_{n-1},y_n)
$$
is an isomorphism, we have 
$$
w\in\la x\cdot y_n\ra =\la y_0y_1\cdots y_n\ra
$$
and hence the result.
\end{proof}

Let $K_n$ be the pure braided group, that is, $K_n$ is that kernel of the
canonical epimorphism $B_n\to\Sigma_n$. Let $B_{n+1}$ be the subgroup of 
$B_{n+2}$ generated by $\sigma_j$ for 
$j\geq0$. Recall that 
$$
\pi_{n+1}(F(S^1))=\calZ(F(S^1))_{n+1}/\calB(F(S^1))_{n+1}.
$$
Consider the actions of two braided groups $B_{n+2}$ and $B_{n+1}$ on
$F(S^1)_{n+1}/\calB(F(S^1))_{n+1}$. We have 
\begin{thm}\label{theorem2.11}
If $n\geq2$, then in $F(S^1)_{n+1}/\calB(F(S^1))_{n+1}$,
\begin{enumerate}
\item[1)] the fixed set of the pure braided group $K_{n+1}$-action
is 
$$
\Z\oplus \pi_{n+1}(F(S^1));
$$
\item[2)] the  fixed set of $K_{n+2}$-action is 
$$\pi_{n+1}(F(S^1)).$$
\end{enumerate}
\end{thm}
\begin{proof}
We show that the homotopy fixed set of $K_{n+1}$ on $F(S^1)_{n+1}$
is generated by $y_{-1}$ and $\calZ(F(S^1))_{n+1}$. Assertions~1 and~2 will
follow from this statement. 
Let $w$ be a homotopy fixed point of $K_{n+1}$-action on $F(S^1)_{n+1}$. 
Since 
$$
\sigma_k^2d_0=d_0\sigma_{k+1}^2
$$
for $k\geq-1$, we have
$$
\sigma_k^2d_0(w)=d_0(w)
$$
for each $k\geq-1$. By Lemma~\ref{lemma2.10}, we have
$$
d_0(w)=1.
$$
Now for each 
$1\leq j\leq n+1$, we have
$$
\sigma_k^2d_j=\left\{
\begin{array}{cccc}
d_j\sigma_{k+1}^2&\quad{\rm if}\quad& j\leq k+1;\\
d_j\circ \sigma_{j-1}^{-1}\circ\sigma_{j-2}^2\circ\sigma_{j-1}&
\quad{\rm if}\quad& j=k+2;\\
d_j\sigma_k^2&\quad{\rm if}\quad& j>k+2.\\
\end{array}
\right.
$$
By Lemma~\ref{lemma2.10}, there exists integers $k_1,k_2,\cdots,k_{n+1}$
such that
$$
d_j(w)=y_{-1}^{k_j}
$$
for $1\leq j\leq n+1$. Since $d_ky_{-1}=y_{-1}$ for $k>0$, we have
$$
y_{-1}^{k_j}=d_j(y_{-1}^{k_j})=d_jd_jw=d_jd_{j+1}w=d_j(y_{-1}^{k_{j+1}})=
y_{-1}^{k_{j+1}}
$$
for $1\leq j\leq n$ and so 
$$
k_1=k_2=\cdots =k_{n+1}.
$$
Let $w'=y_{-1}^{-k_1}w$. Then
$$
d_j(w')=1
$$
for each $0\leq j\leq n+1$ and $w'\in \calZ_(F(S^1))_{n+1}$. This shows that
$w$ lies in the subgroup generated by $y_{-1}$ and cycles 
$\calZ(F(S^1))_{n+1}$ and hence the result.
\end{proof}

\noindent{\bf Note:}
In $F(S^1)_{n+1}/\calB_{n+1}$,
since any element in the homotopy group is a 
homotopy fixed point of $\tilde B_{n+2}$, the 
fixed set of $\tilde B_{n+1}$ is $\Z\oplus \pi_{n+1}(F(S^1))$ and the 
 fixed set of $\tilde B_{n+2}$ is $\pi_{n+1}(F(S^1))$.

\section{Braided Representation of $F(S^1)$}\label{section3}
\subsection{A Representation of $F(S^1)$}
In this subsection, the ground ring $R$ is $\Z$ or $\Z/p$.
Let $X$ be a pointed set. Let $A(X)$ be the algebra of non-commutative
 formal power series 
with variables in any finite subset of $X$ over $R$ modulo the single 
relation that 
$\ast=1$, where $\ast$ is the base point of $X$. Let $X$ be a pointed 
simplicial set. The simplicial algebra $A(X)$ is defined by applying the 
functor $A$ to $X$. Let $X=S^1$ be the simplicial circle. By using the methods
in~\cite{Wu1}, we have
\begin{enumerate}
\item[1)] There is a choice of generators in $A(S^1)_{n+1}$ such that
 $$A(S^1)_{n+1}=A(x_0,x_1,\cdots,x_n)$$ 
is the associated algebra of the 
non-commutative formal
power series in variables $x_0,x_1,\cdots,x_n$ over $R$.
\item[2)] The simplicial 
structure on $A(S^1)$ is given by
$$
d_jx_k=\left\{\begin{array}{cccccc}
x_{k-1}&\mbox{ for}& \mbox{ $j\leq k,$}\\
0&\mbox{ for}&\mbox{ $j=k+1,$}\\
x_k&\mbox{ for} &\mbox{ $j>k+1,$}\\
\end{array}
\right.
$$
and
$$
s_jx_k=\left\{\begin{array}{cccccc}
x_{k+1}&\mbox{ for}&\mbox{ $j\leq k,$}\\
x_k+x_{k+1}+x_kx_{k+1}&\mbox{ for}&\mbox{ $j=k+1,$} \\
x_k&\mbox{ for}&\mbox{ $j>k+1$,}\\
\end{array}
\right.
$$
for $0\leq j\leq n+1$, where
$$
x_{-1}=(1+x_{n-1})^{-1}(1+x_{n-2})^{-1}\cdots (1+x_0)^{-1}-1
$$
in $A(S^1)_n$.
\end{enumerate} 

Let the Braided group $B_{n+1}$ act on $A(S^1)_{n+1}$ as follows.
\begin{enumerate}
\item[] For each $0\leq k\leq n-1$, $\sigma_k\colon A(S^1)_{n+1}\to
A(S^1)_{n+1}$ is an automorphism of algebras with
$$
\sigma_k(x_j)=\left\{
\begin{array}{ccc}
x_{k+1}&\quad {\rm if}&\quad j=k\\
(1+x_{k+1})^{-1}(1+x_k)(1+x_{k+1})-1&\quad {\rm if}&\quad j=k+1\\
x_j & &{\rm otherwise}\\
\end{array}
\right.
$$
for $0\leq k\leq n-1$. 
\end{enumerate}
 Let $\sigma_{-1}$ be an automorphism of 
$A(S^1)_{n+1}$ defined by
$$
\sigma_{-1}(x_0)=(1+x_0)^{-1}(1+x_{-1})(1+x_0)-1=
(1+x_0)^{-1}(1+x_n)^{-1}\cdots (1+x_1)^{-1}-1
$$
and $\sigma(x_j)=x_j$  for $j>0$.
The subgroup of the automorphism group of $A(S^1)_{n+1}$ generated by
$\sigma_j$ for $-1\leq j\leq n-1$ is the braid group $B_{n+2}$. Let
$$
e\colon F(S^1)\to A(S^1)
$$
be the canonical representation, that is
$$
e(y_j)=1+x_j
$$
for each $j$.

\begin{prop}\label{proposition3.1}
The simplicial algebra $A(S^1)$ is a braided simplicial algebra, that is,
the braided action satisfies the identities~\ref{equation1} 
and~\ref{equation2} in Proposition~\ref{proposition2.1}. 
Furthermore, the representation $e\colon F(S^1)\to
A(S^1)$ is a braided representation, that is $e$ commutes with the braid group
action and the simplicial structure.
\end{prop}

The proof is straight forward.

Since the representation $e\colon F(S^1)\to A(S^1)$ is faithful~\cite{MKS}, 
we have

\begin{prop}\label{proposition3.2}
Let $w\in F(S^1)_{n+1}$. Then
\begin{enumerate}
\item[1)] $w\in\calZ(F(S^1))_{n+1}$ if and only if 
$e(w)-1\in\calZ(A(S^1))_{n+1}$;
\item[2)] $w\in NF(S^1)_{n+1}$ if and only if $e(w)-1\in NA(S^1)_{n+1}$.
\end{enumerate}
\end{prop}

Let $w=x_{i_1}^{n_1}x_{i_2}^{n_2}\cdots x_{i_t}^{n_t}$ be a monomial in
$A(x_0,\cdots, x_n)$. We call $w$ is {\it non-degenerate} if the set
$$
\{i_1,\cdots,i_t\}=\{0,\cdots,n\},
$$
that is each letter $x_j$ with $0\leq j\leq n$ appears in $w$ at least once.

\begin{thm}\label{theorem3.3}
Let $f$ be a series in $A(S^1)_{n+1}$. Then $f\in NA(S^1)_{n+1}$ if and only
if $f$ is a formal series of non-degenerate monomials.
\end{thm}

The assertion follows from the following lemma.

\begin{lem}\label{lemma3.4}
Let $f$ be a series in $A(S^1)_{n+1}$. Then 
$$
f\in \bigcap_{j=1}^{i+1}\Ker(d_j)
$$
if and only if $f$ is a linear summation of monomials in which each $x_j$ 
appears at least once for $0\leq j\leq i$.
\end{lem}
\begin{proof}
The proof is given by induction on $i$. Let $i=0$. Since $d_1$ is the 
projection, $\Ker (d_0)$ is a two sided ideal generated by $x_0$. Suppose that
the assertion holds for $i-1$ with $i>0$. Then there is a decomposition
$$
\bigcap_{j=1}^i\Ker d_j=C\oplus D,
$$
where $C$ is the set of series of monomials in which each $x_j$ appears at 
least once for $0\leq j\leq i$. Since $d_{i+1}$ is the projection, we find 
that $d_{i+1}|_C=0$ and $d_{i+1}|_D$ is a monomorphism from $D$  to $A(S^1)_n$.
This shows that 
$$
\bigcap_{j=1}^{i+1}\Ker d_j=C
$$
and hence the result.
\end{proof}

\subsection{Formal Steenrod Operations and Higher Differentials}
In this subsection, the ground ring $R$ is $\Z$ or $\Z/p$.
Let $A(x_0,\cdots,x_n)_i$ be the sub $R$-module of $A(x_0,\cdots,x_n)$
generated by monomial of degree $i$. Let $f\colon A(x_0,\cdots,x_n)\to
 A(x_0,\cdots,x_{n-1})$ be an $R$-linear map. We call $f$ is a homogenous
map of degree $t$ if
$$
f(A(x_0,\cdots,x_n)_i)\subseteq A(x_0,\cdots,x_{n-1})_{i+t}
$$
for each $q$. Let $V$ be the free $R$-module generated by $x_0,\cdots,x_n$.

\begin{lem}\label{lemma3.5}
Let $f_t\colon V\to A(x_0,\cdots x_{n-1})$ be a sequence of $R$-linear
maps with $t\geq0$ such that 
$$
f_t(V)\subseteq A(x_0,\cdots,x_{n-1})_t
$$
for each $t$. Then there exists a unique sequence of homogeneous 
maps $P_f^t\colon A(x_0,\cdots,x_n)\to A(x_0,\cdots,x_{n-1})$ such that
\begin{enumerate}
\item[1)] $P_f^t|_V=f_t$ for each $t\geq0$;
\item[2)] The anti-Cartan formula
$$
P_f^t(xy)=\sum_{i+j=t} P_f^i(y)P_f^j(x)
$$
hold for any $x,y\in A(x_0,\cdots,x_n)$.
\end{enumerate}
\end{lem}

The proof is straight forward.

Let $\chi\colon A(x_0,\cdots,x_n)$ be the convolution, that is, $\chi$ is
the anti-automorphism with $\chi(x_i)=-x_i$ for each $i$. Let 
$$
\bar d_0=d_0\circ\chi=\chi\circ d_0.
$$
Then we have $\bar d_0(x_j)=x_{j-1}$ for $j>0$ and
$$
\bar d_0(x_0)=(1+x_0)(1+x_1)\cdots (1+x_{n-1})-1.
$$
Let 
$$
\Delta_{s-1}=\sum_{0\leq l_1<l_2<\cdots<l_s\leq n-1}x_{l_1}x_{l_2}\cdots 
x_{l_s}.
$$
for $1\leq s\leq n$. For $0\leq s\leq n-1$,  
let $\delta_s\colon V\to A(x_0,\cdots,x_{n-1})$
be the $R$-linear map defined by 
$$
\delta_i(x_j)=\left\{
\begin{array}{ccc} 
x_{j-1}&\quad{\rm for\quad}&  i=0,j>0,\\
0&\quad{\rm for\quad}&i>0,j>0\\
\Delta_i&\quad{\rm for}\quad& j=0\\
\end{array}
\right.
$$

Let
$$
\delta=(\delta_0,\delta_1,\cdots,\delta_{n-1},0,\cdots,0).
$$
\begin{prop}\label{proposition3.6}
The map $\bar d_0\colon A(x_0,\cdots,x_n)\to A(x_0,\cdots,x_{n-1})$ is 
decomposed as
$$
\bar d_0=\sum_{i=0}^{\infty}P_{\delta}^i.
$$
\end{prop}

\begin{proof}
By the anti-Cartan formula, the map
$$
\sum_{i=0}^{\infty}P_{\delta}^i
$$
is a well-defined anti-homomorphism of algebras. Since
$$
\bar d_0(x_j)=\sum_{i=0}^{\infty}P_{\delta}^i(x_j)
$$
for each $0\leq j\leq n$, the assertion follows.
\end{proof}

Let $q_i\colon A(x_0,\cdots,x_n)\to A(x_0,\cdots,x_n)$ be the composite
$$
A(x_0,\cdots,x_n)\rTo^{\rm proj.} A(x_0,\cdots,x_n)_i\rInto A(x_0,\cdots,x_n)
$$
for $0\leq i<\infty$. Let 
$$
\Gamma_tF(S^1)_{n+1}=\{w\in F(S^1)_{n+1}|q_i(w)=0\quad{\rm for}\quad 0<i<t\}
$$
for $1\leq t<\infty$. It is well-known~\cite{MKS} 
that $\{\Gamma_t\}$ is the descending
central series if $R=\Z$ and descending $p$-central series if $R=\Z/p$.
Let
$$
L(S^1)=\bigoplus_{t=1}^{\infty}\Gamma_t(F(S^1))/\Gamma_{t+1}(F(S^1)).
$$
If $R=\Z$, then $L(S^1)$ is the free simplicial Lie algebra over $\Z$ 
generated by $S^1$. If $R=\Z/p$, then 
$L(S^1)$ is the free restricted simplicial Lie algebra over $\Z/p$ generated
by $S^1$. The simplicial
structure  $L(S^1)$ is as follows.
\begin{enumerate}
\item[1)]  $L(S^1)_{n+1}=L(x_0,\cdots,x_n)$;
\item[2)] $d_j(x_k)=x_k$ for $k<j-1$, $d_j(x_{j-1})=0$ and $d_j(x_k)=x_{k-1}$
for $k>j-1$, where $x_{-1}=-(x_0+x_1+\cdots+x_{n-1})\in L^u(S^1)_n$ or
$L(S^1)_n$. 
\item[3)] $s_jx_k=x_{k}$ for $k<j-1$, $s_jx_{j-1}=x_{j-1}+x_j$ and
$s_jx_k=x_{k+1}$ for $k>j-1$.
\end{enumerate}
A Lie monomial $w=[[x_{i_1},x_{i_2},\cdots]^{p^s}$ in $L^u(S^1)_{n+1}$ or 
$L(S^1)_{n+1}$
is called {\it non-degenerate} if the set
$$
\{i_1,\cdots,i_t\}=\{0,1,\cdots,n\}.
$$
By Theorem~\ref{theorem3.3}, we have

\begin{cor}\label{corollary3.7}
The Moore chain complexes  $NL(S^1)_{n+1}$ is the 
submodules of  $L(S^1)_{n+1}$ spanned by non-degenerate
Lie monomials, respectively.
\end{cor}

The spectral sequence, which is denoted by $\{E^r(F(S^1)\}$, 
induced by the descending $p$-central 
(integral descending central) series
of $F(S^1)$ is called (integral) {\it Adams spectral sequence} of $F(S^1)$. 
By Curtis Theorem , this spectral sequence is convergent to
$\pi_*(F(S^1);p)=\pi_{*+1}(S^2;p)$ ( or $\pi_{*+1}(S^2)$ if we use integral
descending central series). A description of higher differentials in the Adams
spectral sequence is as follows.

Let $w\in F(S^1)_{n+1}$. Then 
$$
e(w)=1+\sum_{i=1}^{\infty}q_i(w)
$$
in $A(x_0,\cdots,x_n)$. We simply write $(w)_i$ for $q_i(w)$. Let 
$$
\alpha\in E^1_{t,*}(F(S^1))=\pi_*(L_t(S^1))
$$
and let $z^{\alpha}$ be a 
cycle in the simplicial group $L(S^1)$ such that the homotopy
class of $z^{\alpha}$ is $\alpha$, that is $z^{\alpha}$ is a 
cycle representative of $\alpha$.
Since the map
$$
g\colon\Gamma_t(F(S^1))\to L_t(S^1)=\Gamma_t(F(S^1))/\Gamma_{t+1}(F(S^1)),
$$
is a simplicial epimorphism, there is an element 
$$
w^{\alpha}\in N(\Gamma_t(F(S^1)))=N(F(S^1))\cap F(S^1)
$$ such that 
$$
g(w^{\alpha})=z^{\alpha}.
$$
The element $w^{\alpha}$ is called a {\it Moore representative} of $\alpha$.

\begin{thm}\label{theorem3.8}
Let $\alpha\in E^1_{t,*}(F(S^1_{n+1}))$ and let $w_{\alpha}$ be a Moore 
representative in $F(S^1)_{n+1}$. Let $1\leq r\leq \infty$. 
Then $d^j(\alpha)=0$ for $j<r$ if and only if the following linear equations
holds in $A(x_0,\cdots,x_{n-1})$:
$$
\sum_{i=0}^jP_{\delta}^i w^{\alpha}_{t+j-i}=0.
$$
for $0\leq j<r$. Furthermore, if $r<\infty$, then 
$$
\sum_{i=0}^rP_{\delta}^i w^{\alpha}_{t+r}\in L_{t+r}(S^1)_n\subseteq 
A(x_0,\cdots,x_{n-1}),
$$
which is a cycle representative of $-d^r(\alpha)$.
\end{thm}
\begin{proof}
Since $d^r(\alpha)=0$ if and only if 
$$
d_0(w^{\alpha})\in\Gamma_{t+r}F(S^1)_n,
$$
the assertion follows from Proposition~\ref{proposition3.6}.
\end{proof}

Let $w$ be a word in $F(S^1)_{n+1}$. We call $w$ is a {\it 
basic non-degenerate 
commutator} of weight $s=l(w)$ if $w$ can be written down as a commutator
$$
[\cdots [y_{i_1},y_{i_2}],\cdots],y_{i_s}]
$$
 such that the set
$$
\{i_1,\cdots,i_s\}=\{0,\cdots,n\},
$$
that is each generator $y_j$ occurs in the commutator $w$ at least once. 
Let $N_tF(S^1)_{n+1}$ be the subgroup of $F(S^1)_{n+1}$ generated by 
basic nondegenerate commutators $w$ with $l(w)\geq t$ and let 
$N_t^{(p)}F(S^1)_{n+1}$ be the subgroup of $F(S^1)_{n+1}$ generated by 
$
w^{p^r},
$
where $w$ runs over all basic 
non-degenerated commutators with $l(w)\cdot p^r\geq t$. Note that
$$
N_tF(S^1)_{n+1}\subseteq N(F(S^1))_{n+1}\cap \Gamma_t(F(S^1)_{n+1})
$$
if $R=\Z$ and
$$
N_t^{(p)}F(S^1)_{n+1}\subseteq N(F(S^1))_{n+1}\cap\Gamma_t(F(S^1)_{n+1})
$$
if $R=\Z/p$.

By Corollary~\ref{corollary3.7}, we have

\begin{prop}\label{proposition3.9}
Let $r$ be any non-negative integer. If $R=\Z$, there is an isomorphism
$$
N_tF(S^1)_{n+1}/(N_tF(S^1)_{n+1}\cap \Gamma_{t+r}F(S^1)_{n+1})
$$
$$
\cong
N(F(S^1))_{n+1}\cap\Gamma_t(F(S^1))_{n+1}/(N(F(S^1))_{n+1}\cap
\Gamma_{t+r}(F(S^1))_{n+1}).
$$
If $R=\Z/p$, then there is an isomorphism
$$
N_t^{(p)}F(S^1)_{n+1}/(N_t^{(p)}F(S^1)_{n+1}\cap \Gamma_{t+r}F(S^1)_{n+1})
$$
$$
\cong
N(F(S^1))_{n+1}\cap\Gamma_t(F(S^1))_{n+1}/(N(F(S^1))_{n+1}\cap
\Gamma_{t+r}(F(S^1))_{n+1}).
$$
\end{prop}

\begin{rem}{\rm 
 By Theorem~\ref{theorem3.8} and 
Proposition~\ref{proposition3.9}, $\pi_*(L_t(S^1))$ are represented by those
words $w$ in $N_tF(S^1)$ (or $N_t^{(p)}F(S^1)$) with $P^0_{\delta}w_t=0$ and
the higher differentials in the Adams spectral sequence are related to higher
formal Steenrod operations on $N_tF(S^1)$. 
}
\end{rem}

\subsection{$E^1$-terms of The Integral Adams Spectral Sequence}
In this subsection, the ground ring $R$ is a subring of $\Q$. 
Let $L$ be the free functor from free $R$-modules to Lie algebras. Let $X$ be
a pointed simplicial set. 
Let $\bar R(X)=R(X)/R(\ast)$ be the reduced free simplicial $R$-module
generated by $X$. In particular, $\bar R(S^n)=K(R,n)$. Let 
$L(X)=L(\bar R(X))$ for any pointed simplicial set $X$.
Let $V$ be a free $R$-module and let $L'$ be the kernel of the 
abelianizer
$$
L(V)\to V.
$$
Then $L'$ is  functor from free $R$-modules to graded free $R$-modules.
Let $Q_n(L'(V))$ be the set of indecomposable elements of degree $n$
of  $L'(V)$. Let $S_n(V)$ be the set of monomials of degree $n$ in the 
polynomial algebra $S(V)$. Note that $Q_1(L'(V))=0$.

\begin{lem}\label{lemma3.10}
For each $n\geq 2$, there is a functorial short exact sequence
$$
0\to Q_n(L'(V))\to S_{n-1}(V)\otimes V\rTo^{\rm mult.} S_n(V)\to0.
$$
\end{lem}
\begin{proof}
Let $K(V)$ be the kernel of $S_{n-1}(V)\otimes V\rTo^{\rm mult}S_n(V)$ and 
let
 $\phi\colon T_n(V)\to S_n(V)$ be the composite
$$
T_n(V)=T_{n-1}(V)\otimes V\rTo^{\rm proj.} S_{n-1}(V)\otimes 
V\rTo^{\rm mult.} S_n(V).
$$
Then $\phi|_{L'_n(V)}\colon L'_n(V)\to S_n(V)$ is zero which gives a functorial
map
$$
\tilde \phi\colon L'_n(V)\to K(V).
$$
It is a routine work to show that the map $\tilde\phi$ is onto and
$\tilde \phi|_{D_n(L'(V))}\colon D_n(L'(V))\to K(V)$ is zero and so 
$\tilde \phi$ factors through $Q_n(L'(V))$, where $D(L'(V))$ is the set of 
decomposable elements of $L'(V)$. Thus $K(V)$ is a functorial quotient of
$Q_n(L'(V))$. By checking Poincare series, the quotient map
$$
Q_n(L'(V))\to K(V)
$$
is an isomorphism and hence the result.
\end{proof}

\noindent{\bf Note:} There is no functorial cross-section from $Q(L')$ to
$L'$ (See~\cite{SW}).

\begin{prop}\label{proposition3.11}
There is a homotopy equivalence
$$
L(S^1)\simeq \bar R(S^1)\oplus L(S^2).
$$
\end{prop}
\begin{proof}
Since $K(\Z,1)\simeq S^1$, we have 
$$
Q(L')\simeq\bar R(S^2).
$$
Let $\phi\colon L(L_2(S^1))\to L'(S^1)$ be the inclusion. Then $\phi$ is a 
homotopy equivalence by checking the spectral sequence induced by Lie 
filtrations. The assertion follows.
\end{proof}

Let $V$ be a free simplicial $R$-module and let $C$ be a (pointed)
 simplicial coalgebra. Let $f,g\colon
C\to T(V)$ be pointed simplicial coalgebra maps. We call $f$ is 
{\it coalgebra homotopic}
to $g$ if there is a pointed homotopy $F_t\colon C\to T(V)$ such that $F_0=f$, 
$F_1=g$ and $F_t$ is a coalgebra map for each $t$. Let
$$
[T(V),T(V)]^{\rm coalg}
$$
be the set of coalgebra homotopy classes. If $f\colon T(V)\to T(V)$ is a 
simplicial coalgebra map, then we have the restriction
$$
f|_{L(V)}\colon L(V)\to L(V).
$$
If  $f$ is coalgebra homotopic to $g$, then 
$$
f|_{L(V)}\simeq g|_{L(V)}.
$$
This defines a map
$$
\theta\colon [T(V),T(V)]^{\rm coalg}\to [L(V),L(V)].
$$

\begin{thm}\label{theorem3.12}
Suppose that $R=\Z_{(p)}$. Then
the homotopy groups $\pi_*(L(S^{2n}))$ has the following exponents.
\begin{enumerate}
\item[1)] If $n=1,3$, then
$$
p\cdot \pi_*(L(S^{2n}))=0
$$
for $*>2n$ and any prime $p$;
\item[2)] If $p=2$, then
$$
4\cdot \pi_*(L_t(S^{2n}))=0
$$
for $t>2$ and any $n$.
\end{enumerate}
\end{thm}
\begin{proof}
Consider the Cohen representation
$$
\theta\colon H_{\infty}\to [T(S^{2n}),T(S^{2n})]^{\rm coalg}\to 
[L(S^{2n}),L(S^{2n})],
$$
where $H_{\infty}$ is the Cohen group. 
First we assume that $n=1,3$.  Then the Samelson product
$$
S^{2n}\wedge S^{2n}\to F(S^{2n})
$$
is null homotopic. Since $\pi_*(\Gamma_2(F(S^{2n}))\to \pi_*(F(S^{2n}))$ is
a monomorphism, the Samelson product
$$
S^{2n}\wedge S^{2n}\to\Gamma_2(F(S^{2n}))
$$
is null homotopic and so the composite
$$
S^{2n}\wedge S^{2n}\to \Gamma_2(F(S^{2n}))\to\Gamma_2(F(S^{2n})/
\Gamma_3(F(S^{2n}))
$$
is null homotopic. Let $J_tT(S^{2n})=\Z_{(p)}(J_t(S^{2n}))$. Then there is
 a monomorphism
$$
[J_tT(S^{2n}),T(S^{2n})]^{\rm coalg}\to [\Z_{(p)}((S^{2n})^{\times t}), 
T(S^{2n})]^{\rm coalg}
$$
for any $t$. Since the Samelson product is trivial, the group
$$
[\Z_{(p)}((S^{2n})^{\times t}),T(S^{2n})]^{\rm coalg}
$$
is abelian. Let $\id^{\ast p}\colon T(V)\to T(V)$ be the $p$-fold  
convolution product
of the identity and let $T(p)\colon T(V)\to T(V)$ be the morphism of Hopf 
algebras induced by  the map $p\colon V\to V,\quad, x\to px$. Then
$\id^{\ast p}|_{J_tT(S^{2n})}$ is coalgebra homotopic to 
$T(p)|_{J_tT(S^{2n})}$ in
$$
[\Z_{(p)}((S^{2n})^{\times t}),T(S^{2n})]^{\rm coalg}
$$
for each $t$ and so 
$$
\id^{\ast p}|_{J_tT(S^{2n})}\simeq T(p)_{J_tT(S^{2n})}
$$
for each $t$. This $\id^{\ast p}$ is coalgebra homotopic to $T(p)$.
Let $x\in L_t(S^{2n})$. Then 
$$
\id^{\ast p}(x)=px\quad T(p)(x)=p^tx.
$$
Thus
$$
p^t-p=p(p^{t-1}-1)=0
$$
in $[L_t(S^{2n}),L_t(S^{2n})]$
and so
$$
p\cdot \pi_*(L_t(S^{2n}))=0
$$
for $t>1$. Since $L_1(S^{2n})=K(Z_{(p)},2n)$, we proves assertion~1.

The proof of assertion~2 is similar to assertion~1, where one needs the fact
that the Samelson product
$$
S^{2n}\wedge S^{2n}\to F(S^{2n})
$$
is of order $2$ up to homotopy and the higher Samelson products are null 
homotopic localized at $2$.
\end{proof}

\begin{rem}{\rm  $\pi_*(L(S^n)\otimes \Z/p)$ is known as a 
specific module over the $\Lambda$-algebra~\cite{Bousfield,Curtis,Wellington}.
 By considering the Bockstein 
spectral sequence for $\pi_*(L(S^n))$, this theorem shows that 
$\pi_*(L(S^{2n})$ is the kernel of the Bockstein on $\pi_*(L(S^{2n})\otimes
\Z/p)$ when $n=1,3$.  }
\end{rem}

\begin{thm}\label{theorem3.14}
Let $R=\Z_{(p)}$.
If $t$ is not a power of $p$, then
$$
\pi_*(L_t(S^{2n}))=0.
$$
\end{thm}
\begin{proof}
By a result in~\cite{SW}, there exist functors $A^{\min}$
and $Q^{\max}_t$ for $t\geq2$ such that
\begin{enumerate}
\item[1)] $A^{\min}$ is a (smallest) coalgebra retract of $T$ with
$V\subseteq A^{\min}(V)$ for each $V$;
\item[2)] $Q^{\max}_t$ is a subfunctor of $L_t$ and $Q^{\max}_t$ is a 
retract of $T_t$;
\item[3)] There is a functorial coalgebra decomposition
$$
T(V)\cong T(\bigoplus_{t=2}^{\infty}Q^{\max}_t(V))\otimes A^{\min}(V).
$$
\end{enumerate}
Let $L^{\min}(V)$ be the primitives of $A^{\min}(V)$. Then there is a 
functorial decomposition
$$
L(V)\cong L(\bigoplus_{t=2}^{\max}Q^{\max}(V))\oplus L^{\min}(V).
$$
Since $Q^{\max}_t$ is a retract of $T_t$, $Q^{\max}(S^{2n})$ is a retract
of $T_t(S^{2n})$ and so $Q^{\max}_t(S^{2n})$ is either contractible or 
homotopic to $K(\Z_{(p)}, 2tn)$. Now $Q^{\max}_t(S^{2n})$ is contractible
because $Q^{\max}_t$ is subfunctor of $L_t(V)$. It follows that 
$$
L(\bigoplus_{t=2}^{\infty}Q^{\max}_t(S^{2n}))
$$
is contractible. The assertion follows from a result in~\cite{SW} that
$L_t$ is a functorial retract of 
$$
L(\bigoplus_{t=2}^{\infty}Q^{\max}_t(S^{2n}))
$$
if $t$ is not a power of $p$.
\end{proof}

By Theorems~\ref{theorem3.12} and~\ref{theorem3.14}, we have

\begin{thm}
Let $\{E^r\}_{r\geq1}$ 
be the integral Adams spectral sequence of $F(S^1)$. Then
\begin{enumerate}
\item[1)] $E^r_{t,*}=0$ unless $t=2p^s$ some prime $p$ and some 
non-negative integer $s$.
\item[2)] $p\cdot E^r_{2p^s,*}=0$ for any prime $p$ and any integer $s>0$.
\item[3)] Let $\alpha\in E^1_{2p^s,*}$ with $s>0$. Then the only  
differentials $d^r$, which are possibly non-trivial on $\alpha$, are
$$
d^{2p^t-2p^s}
$$
with $t>s$.
\end{enumerate}
\end{thm}

\end{document}